\pdfoutput=1
\documentclass[11pt]{amsart}
\usepackage[T1]{fontenc}
\usepackage[margin=.95in]{geometry}
\usepackage{amsmath,amssymb,amsthm}
\usepackage{etoolbox}
\makeatletter
\patchcmd{\@settitle}{\uppercasenonmath\@title}{\Large}{}{}
\patchcmd{\@setauthors}{\MakeUppercase{\authors}}{\authors}{}{}
\patchcmd{\maketitle}{\uppercasenonmath\shorttitle}{}{}{}
\patchcmd{\maketitle}{\@nx\MakeUppercase{\the\toks@}}{\the\toks@}{}{}
\makeatother
\usepackage{hyperref}
\usepackage{esint}

\newtheorem{theorem}{Theorem}[section]
\newtheorem{proposition}[theorem]{Proposition}
\newtheorem{lemma}[theorem]{Lemma}
\newtheorem{corollary}[theorem]{Corollary}
\theoremstyle{remark}\newtheorem{remark}[theorem]{Remark}
\newtheorem*{remark*}{Remark}
\newtheorem*{note*}{Note}
\newcommand{\E}{\mathcal E}
\newcommand{\R}{\mathbb R}
\newcommand{\Ric}{\operatorname{Ric}}
\newcommand{\Vol}{\operatorname{Vol}}
\newcommand{\im}{\operatorname{im}}
\newcommand{\II}{\mathrm{II}}
\newcommand{\dvol}{\,d\mathrm{vol}}

\title{Energy and topology of Ricci-flat $4$-manifolds}
\author{Tristan Ozuch}
\address{MIT, Dept. of Math., 77 Massachusetts Avenue, Cambridge, MA 02139-4307}
\email{ozuch@mit.edu}
\date{}
\begin{document}

\begin{abstract}
    We show that there is $C>0$ such that any complete Ricci-flat $4$-manifold $(M,g)$ satisfies
$$
b_2(M)-C\leqslant\frac1{8\pi^2}\int_M|R_g|_g^2\,dv_g\leqslant b_2(M)+C,
$$
where $R$ is the curvature operator. No assumption is made on volume growth, model at infinity or curvature decay. It follows that uncountably many exotic $\R^4$ admit no complete Ricci-flat metric, and that complete noncompact hyperkähler metrics with finite $b_2$ have ends classified by Sun-Zhang. This also provides information about the thin part of immortal four-dimensional Ricci flows in Bamler's decomposition. The proof refines the chopping method of Cheeger-Gromov and Cheeger-Tian by bounding the topology of boundaries of regions exhausting the manifold.
\end{abstract}

\maketitle

\section*{Introduction}
\renewcommand{\thetheorem}{\Alph{theorem}}

On a closed Riemannian four-manifold $(M,g)$, denote the curvature operator $R\colon\Lambda^2\to\Lambda^2$ and define the \textit{energy} $\E(U):=(8\pi^2)^{-1}\int_U |R|^2\dvol$. The Chern-Gauss-Bonnet formula gives an identity between the energy and the topology of $M^4$
\[
 \E(M)=\chi(M)+\frac1{8\pi^2}\int_M\bigl|\Ric-\tfrac14\operatorname{Scal}\,g\bigr|^2\dvol .
\]
In particular, we have $\E(M)\geqslant\chi(M)$, with equality precisely for Einstein metrics which minimize $\mathcal{E}$. This is one of the motivations for Hawking's name \emph{gravitational instanton} \cite{Haw77}. On noncompact Ricci-flat four-manifolds, the definition of a gravitational instanton often assumes $\E(M)<+\infty$.

In this article, we show that having a finite second Betti number is equivalent to this finite energy condition and answer several open questions from the literature, e.g. \cite{LS26GI}. The result is already known for a complete noncompact \textit{noncollapsed} Ricci-flat four-manifold $M$: \cite{CN15} and \cite{BKN89} give an ALE end asymptotic to a flat cone $\R^4/\Gamma$ and in this case $\E(M)=\chi(M)-\frac{1}{|\Gamma|}$ by \cite{Nak90}. 

\begin{theorem}\label{thm:main}
    There exists $C>0$ with the following properties. Let $(M,g)$ be a complete noncompact Ricci-flat $4$-manifold with $b_2(M)<+\infty$ and $p\in M$. The energy $\mathcal{E}(M)$ satisfies $$b_2(M)-C\leqslant\mathcal{E}(M)\leqslant b_2(M)+C.$$
    
    \noindent Additionally, the following holds if $\lim_{r\to+\infty}r^{-4}\Vol B_r(p)=0$.
    \begin{enumerate}
        \item Outside a compact subset of $M$, we have $|R|\leqslant C\,d(\cdot,p)^{-2}$. For every sufficiently large $r$, $B_{2r}(p)\setminus\overline{B_r(p)}$ carries a \emph{nilpotent structure} of positive rank whose orbits have extrinsic diameter $o(r)$ in the sense of \cite{CFG92}.
        \item For every sufficiently large $r$, there is a compact region $B_r(p)\subset Z_r\subset B_{2r}(p)$ with smooth boundary $\partial Z_r$ such that $b_0(\partial Z_r)+b_1(\partial Z_r)\leqslant C $, $b_2(Z_r)\leqslant\E(M)+C$ and  $$\E(M)=\chi(Z_r)\in \mathbb Z.$$
        \item If $M$ is oriented, the eta invariants of the signature operators of $\partial Z_r$ converge to a finite value $ \eta_\infty$. With $\tau(M)$ the signature of the intersection form of $M$, we have the following generalizations of the Hitchin-Thorpe inequality and a bound on $\eta_\infty$:$$3\left|\eta_\infty+\tau(M)\right|\leqslant 2\,\E(M)\leqslant 2\,b_2(M)+C
\qquad\text{and}\qquad
3\,|\eta_\infty|\leqslant 5\,b_2(M)+C.$$
    \end{enumerate}
\end{theorem}

Theorem \ref{thm:main} solves parts of a conjecture of Sun-Zhang and Li-Sun \cite[Conjecture 30]{LS26GI}.

\begin{remark*}
    Infinite energy examples are found in \cite{AKL89}, \cite{Myers,KorotkinNicolai,KRWY} and \cite{Ozuch6}. 
\end{remark*}

This structure first implies that many exotic $4$-manifolds cannot carry Ricci-flat metrics.

\begin{corollary}\label{cor:exotic}
    Complete Ricci-flat metrics cannot exist on uncountably many exotic $\mathbb R^4$. 
\end{corollary}

\begin{remark*}
    In \cite[Theorem 1.3]{CMM24}, Chodosh, Maximo and Mukherjee exclude complete metrics of uniformly positive scalar curvature on uncountably many exotic $\R^4$ with Taylor invariant zero \cite{Tay97}, while Corollary \ref{cor:exotic} excludes Ricci-flat metrics on the smoothings of $\R^4$ with Taylor invariant $\gamma>2C$, which include uncountably many smoothings with $\gamma=+\infty$ by \cite[Example 5.10]{Tay97}. 
\end{remark*}

Theorem \ref{thm:main} lets us additionally remove assumptions from the classification of hyperkähler metrics. Chen and Chen \cite{ChenChen21} treated hyperk\"ahler four-manifolds with faster than quadratic curvature decay while Sun and Zhang \cite{SZ24} obtained the end classification under finite curvature energy. 

\begin{corollary}\label{cor:hyperkahler}
    Every nonflat complete noncompact hyperkähler manifold $(M^4,g)$ with $b_2(M)<+\infty$ is diffeomorphic to the interior of a compact manifold with boundary and its end is ALE, ALF, ALG, ALH, ALG$^*$ or ALH$^*$ in the classification of Sun and Zhang \cite[Theorem 1.2]{SZ24}. 
\end{corollary}

The upper bound in Theorem \ref{thm:main} passes to Ricci-flat metrics bubbling off from other manifolds. 

\begin{corollary}\label{cor:limits}
There exists $C>0$ such that the following holds. Let $(M^4,g,p)$ be a complete Ricci-flat orbifold with finite singular set $S$ and singularity group $\Gamma_q$ at $q\in S$. If closed pointed four-manifolds $(N_j,g_j,p_j)$ converge to $(M,g,p)$ smoothly on compact subsets of $M\setminus S$, then we have
$$\mathcal E(M\setminus S)+\sum_{q\in S}\Big(1-\frac1{|\Gamma_q|}\Big)
\leqslant \liminf_{j\to+\infty} b_2(N_j)+C.$$
\end{corollary}

\noindent This statement applies to the Ricci-flat blow-up limits in the thin part of Bamler's decomposition \cite[Corollary 2.42]{Bam20} of immortal Ricci flows, which is an application of the compactness and partial regularity theory of \cite{Bam20a,Bam23,Bam20}. It would also apply to a conjectural four-dimensional Ricci flow with surgery through which $b_2$ is nonincreasing.  Ricci-flat metrics can already bubble off in the thick part, see \cite{DO24}.

\subsection*{Strategy of proof}

For a compact set $Z\subset M$ with boundary $\partial Z= \Sigma$, we first show that the difference between $b_2(Z)$ and $\mathcal{E}(Z)$ is controlled by the geometric boundary term $J(\partial Z)$ in the Chern-Gauss-Bonnet formula and the Betti numbers $b_i(\partial Z)$ of the boundary $\Sigma^3=\partial Z$. This is proven by analyzing the restrictions of the cohomology of $M$ to $Z$ and $\Sigma$. 

We then localize the chopping construction of Cheeger and Gromov \cite{CG91} at the curvature scale $r_{|R|}$, following Cheeger and Tian \cite[Sections 2-3]{CT06}, which provide a way to choose $Z$ with $J(\partial Z)$ controlled by $\int_\Sigma r_{|R|}^{-3}$. The same bound fails for the boundary Betti numbers (see, e.g. finite tubes in $\mathbb R\times \mathbb T^3$), so we instead control the $b_i(\partial Z)$ by $\int_\Sigma \frac{r_{|R|}}{\Vol(B_{r_{|R|}})}$ through local versions of the controls of harmonic forms and multiplicity of \cite{Li80,LS84}.

The final steps of the proof of Theorem \ref{thm:main} are based on ideas from the proof of the key estimate \cite[Theorem 1.26]{CT06}. The $\varepsilon$-regularity estimate of \cite[Theorem 0.8]{CT06} and a maximal function estimate let us bound the $b_i(\partial Z)$ and $J(\partial Z)$ by $\mathcal{E}(A)^\frac34$ for an annulus $A$ surrounding $\Sigma$. Iterating this sublinear bound proves Theorem \ref{thm:main}. 

\begin{note*}
   We write proofs in the oriented case, but the energy and Betti number estimates also hold when $M$, $N$ and the $N_j$ are nonorientable. All along the article, the constant $C$ depends only on the dimension, here $4$, unless stated otherwise. Its value may change from line to line.
\end{note*}

\subsection*{Acknowledgements}

The author is partially supported by the grant NSF, DMS-2405328. He thanks Richard Bamler for asking him this question and explaining its Ricci flow motivation.

\subsection*{Use of AI tools}

LLMs were used to find references, to proofread  the arguments and grammar of the article, and to help fix small gaps in the literature cited. The application Corollary \ref{cor:exotic} was suggested by an LLM.

\renewcommand{\thetheorem}{\thesection.\arabic{theorem}}
\section{Geometric and topological boundary terms}

Let $Z$ be a compact region with smooth boundary $\Sigma$ in a complete noncompact Ricci-flat four-manifold $(M,g)$. For the second fundamental form $\II$ of $\Sigma$, set
$J(\Sigma):=\int_\Sigma (|R||\II|+|\II|^3)\,dA$.
The Chern-Gauss-Bonnet formula with boundary gives
\begin{equation}\label{eq:gauss-bonnet}
|\E(Z)-\chi(Z)|\leqslant CJ(\Sigma).
\end{equation}

We now show that the difference between $\chi(Z)$ and $b_2(Z)$ is controlled by the topology of its boundary, so we get the following estimate. All cohomology groups are real in this article. 

\begin{lemma}\label{lem:boundary-homology}
Let $Z$ be compact with smooth boundary $\Sigma$ in an oriented four-manifold $N$. Suppose that every component of $Z$ meets $\Sigma$. Then, denoting $H^k_c$ the compactly supported cohomology of degree $k$ as in \cite[Definition 1.3]{Car07}, we have
\[
 \chi(Z)\leqslant \dim H^2_c(N)+b_0(\Sigma)+b_1(\Sigma).
\]
If $Z$ embeds in a connected complete oriented four-manifold $M$ of infinite volume with $\operatorname{Ric}\geqslant0$, then
\[
 b_2(Z)\leqslant\chi(Z)+b_0(\Sigma)+b_1(\Sigma).
\]
\end{lemma}
\begin{proof}
    Let $r_k: H^k(Z)\to H^k(\Sigma)$ be the restriction. Each class in $\ker r_k$ has a closed representative supported in $\mathring{Z}$ by \cite[Section 1.1.1]{Car07}, so in particular $\ker r_k\subset\im(H_c^k(N)\to H^k(Z))$ by extension by zero, hence we have
    \begin{equation}\label{eq:b2Z-control}
         b_k(Z)=\dim\ker r_k+\dim\operatorname{im}r_k\leqslant\dim H^k_c(N)+b_k(\Sigma).
    \end{equation}
Every component of $Z$ meets $\Sigma$ by assumption, so $b_4(Z)=0$ and $b_0(Z)\leqslant b_0(\Sigma)$, and since $\Sigma$ is a closed oriented three-manifold, $\chi(\Sigma)=0$, $b_3(\Sigma)= b_0(\Sigma)$ and $b_2(\Sigma)= b_1(\Sigma)$. The first bound then follows from $\chi(Z)\leqslant b_0(Z)+b_2(Z)$ and \eqref{eq:b2Z-control} with $k=2$.
    
    Suppose now that $M$ is connected, complete with $\operatorname{Ric}\geqslant0$ and infinite volume. By the Bochner argument and \cite[Corollaries 1.6 and 1.12]{Car07}, the maps $H_c^k(M)\to H^k(M)$ vanish for $k\in\{1,3\}$, so $\ker r_k=0$ and $b_1(Z)+b_3(Z)\leqslant b_1(\Sigma)+b_0(\Sigma)$ because $b_k(Z)=\dim\ker r_k+\dim\operatorname{im}r_k$.
\end{proof}

Let now $Z\subset M$ also embed in $N$. Lemma \ref{lem:boundary-homology} together with \eqref{eq:gauss-bonnet} shows that
\begin{equation}\label{eq:reduction}
    \mathcal{E}(Z)\leqslant \dim H^2_c(N)+b_0(\Sigma)+b_1(\Sigma)+CJ(\Sigma) \quad \text{and}\quad b_2(Z)\leqslant \mathcal{E}(Z)+b_0(\Sigma)+b_1(\Sigma)+CJ(\Sigma)
\end{equation}
The rest of the article consists in estimating $b_0(\Sigma)+b_1(\Sigma)+CJ(\Sigma)$. 

\section{Controlling the topology of manifolds at the scale of curvature}

For $a>0$, define the curvature scale $r_{|R|}$ and its truncated version $\ell(x)$ as in \cite{CT06} by
\begin{equation}\label{eq:radius}
r_{|R|}(x):=\sup\bigl\{t>0,\,\sup_{B_t(x)}|R|\leqslant t^{-2}\bigr\} \qquad\text{and}\qquad \ell(x):=\min\{a,r_{|R|}(x)\}.
\end{equation}
The function $\ell$ is $1$-Lipschitz. For $c>0$, put $V_{c\ell}(x):=\Vol B_{c\ell(x)}(x)$. If $|\II|\leqslant C_1\ell^{-1}$ on $\Sigma$, then
\begin{equation}\label{eq:boundary}
J(\Sigma)\leqslant C\int_\Sigma\ell^{-3}\,dA \leqslant C\int_\Sigma \frac{\ell(x)}{V_\ell(x)}dA(x)
\end{equation}
where $C$ depends only on $C_1$. The first bound follows from $|R(x)|\leqslant\ell(x)^{-2}$ and the last one from the Bishop-Gromov inequality. We now control the topology of the boundary. 

\begin{proposition}\label{prop:boundary-topology}
   If $|\II|\leqslant C_1\ell^{-1}$ on $\Sigma$, then there are constants $c\in(0,\frac14)$ and $C>0$, depending only on $C_1$, such that

    \vspace{-0.2cm}
   
    $$b_0(\Sigma)+b_1(\Sigma)\leqslant \int_{\Sigma}\frac{C}{\Vol_\Sigma B^{\Sigma}_{c\ell(x)}(x)} dA(x) \leqslant C\int_{\Sigma}\frac{\ell(x)}{V_\ell(x)}\,dA(x).$$
    
    \vspace{-0.2cm}
\end{proposition}
\begin{proof}
    We have $b_0(\Sigma)\leqslant \int_{\Sigma}\frac{1}{\Vol_\Sigma B^{\Sigma}_{c\ell(x)}(x)} dA(x)$. Fix $x\in\Sigma$ and let $\phi$ be a harmonic $1$-form on $\Sigma$. Since $\ell$ is $1$-Lipschitz, $\ell\geqslant(1-2c)\ell(x)$ on $B^{\Sigma}_{2c\ell(x)}(x)$ and the Gauss-Codazzi equations give $\Ric(g_\Sigma)\geqslant-C\ell(x)^{-2}g_\Sigma$ on $B^{\Sigma}_{2c\ell(x)}(x)$ and we have
 $\Delta_\Sigma|\phi|^2\geqslant-C\ell(x)^{-2}|\phi|^2$ by the Bochner formula. So the mean value inequality of \cite[Corollary 3.6]{Li08} gives
    \vspace{-0.2cm}
    \begin{equation}\label{eq:control-harm-bdry}
        |\phi(x)|^2 \leqslant C\frac{\|\phi\|_{L^2(\Sigma)}^2}{\Vol_\Sigma B^{\Sigma}_{c\ell(x)}(x)}.
    \end{equation}
    
    \vspace{-0.2cm}
    
    \noindent Consider now an $L^2(\Sigma)$-orthonormal basis of harmonic $1$-forms $\phi_i$ on $\Sigma$. Since the span of $\phi_i(x)$ is of rank at most $3$, we can find an orthogonal change of basis in order to have $\phi_i(x)=0$ whenever $i\notin\{1,2,3\}$ while not changing $\sum_i|\phi_i(x)|^2$. In particular, by \eqref{eq:control-harm-bdry}, we find

    \vspace{-0.2cm}
    \begin{equation}\label{eq:sum-forms-control}
        \sum_i|\phi_i(x)|^2 \leqslant C\frac{3}{\Vol_\Sigma B^{\Sigma}_{c\ell(x)}(x)}.
    \end{equation}

    \vspace{-0.2cm}

 \noindent Now $b_1(\Sigma)=\int_\Sigma\sum_i|\phi_i|^2$, by orthonormality of $(\phi_i)_i$, so integrating \eqref{eq:sum-forms-control} over $\Sigma$ gives the first inequality, as in \cite[(2.13)]{Li80}.

    The normal exponential map and standard volume comparison give $V_\ell(x)\leqslant C\ell(x)\Vol_\Sigma B^\Sigma_{c\ell(x)}(x)$ regardless of injectivity for $c=c(C_1)>0$ sufficiently small, which proves the second inequality.
\end{proof}

\section{Choice of the cut}\label{sec:cut}
Let $(M,g)$ be a complete noncompact Ricci-flat four-manifold, $p\in M$ and $0<s\leqslant r$. Put
\[
 \rho:=d(p,\cdot)-r
 \qquad\text{and}\qquad
 A':=B_{r+2s/3}(p)\setminus\overline{B_{r+s/3}(p)}=\{s/3<\rho<2s/3\},
\]
and let $\ell$ be the function defined by \eqref{eq:radius} for a given $a>0$.

\begin{lemma}\label{lem:weighted-cut}
There is $C>0$ such that the following holds. Let $0<a\leqslant s/15$ and let $w\geqslant0$ be a continuous function on $A'$. There is a compact region $Z$ with smooth boundary $\partial Z=\Sigma\subset A'$
such that $B_{r+s/3}(p)\subset Z\subset B_{r+2s/3}(p)$
\begin{equation}\label{eq:weighted-cut}
 |\II|\leqslant C\ell^{-1}\quad\text{on }\Sigma
 \qquad\text{and}\qquad
 \int_\Sigma w\,dA\leqslant\frac Cs\int_{A'}w\dvol.
\end{equation}
\end{lemma}

\begin{proof}
We adapt the chopping construction of \cite[Section 1]{CG91} and \cite[proof of Theorem 3.1]{CT06} using the averaging argument of \cite[proof of Theorem 3.13]{CT06}. Put $I:=[2s/5,3s/5]$ and consider the annular region $U_t:=\{t<\rho<t+\ell\}$ for $t\in I$. Since $\ell\leqslant a\leqslant s/15$, we have $U_t\subset\{2s/5<\rho<2s/3\}\subset A'$. A point $x\in A'$ lies in $U_t$ if and only if $\rho(x)-\ell(x)<t<\rho(x)$, so Tonelli's theorem gives
\[
 \int_I\int_{U_t}\frac w\ell\dvol\,dt
 =\int_{A'}\frac{w(x)}{\ell(x)}\,\bigl|I\cap\bigl(\rho(x)-\ell(x),\rho(x)\bigr)\bigr|\dvol(x)
 \leqslant\int_{A'}w\dvol .
\]
Since $|I|=s/5$, we can find and fix $t_0\in I$ such that
\begin{equation}\label{eq:choice-t}
 \int_{U_{t_0}}\frac w\ell\dvol\leqslant\frac5s\int_{A'}w\dvol .
\end{equation}
For this $t_0$, we introduce a scale-adapted function $h$ which plays the role of $d(\cdot,K)$ in \cite{CG91} in order to apply Yomdin's quantitative Sard theorem uniformly on a covering by balls of radii $\propto\ell$:
\[
 h:=\min\Bigl\{\frac{(\rho-t_0)_+}\ell,\,1\Bigr\}\quad \text{so that}\quad\{0<h<1\}=\{0<\rho-t_0<\ell\}=U_{t_0}.
\]
Since $\rho$ and $\ell$ are $1$-Lipschitz, $|dh|\leqslant2\ell^{-1}$ almost everywhere, so $h$ is $C\ell(q)^{-1}$-Lipschitz on each $B_{\ell(q)/2}(q)$. The smoothing and induction on a covering by balls $B_{\ell(q)}(q)$ in \cite[proof of Theorem 3.1]{CT06}, applied to $h$, give $F\in C^\infty(M)$ and $\delta\in(0,\frac12)$ depending only on the dimension such that 
\begin{equation}\label{eq:one-sided}
 F\leqslant h\leqslant F+\tfrac12\quad\text{on }M,
\end{equation}
\begin{equation}\label{eq:band}
 C^{-1}\ell^{-1}\leqslant|\nabla F|\leqslant C\ell^{-1}
 \qquad\text{and}\qquad
 |\operatorname{Hess}F|\leqslant C\ell^{-2}
 \quad\text{on }\{0<F<\delta\}.
\end{equation}
\begin{remark}
    We use the estimate \eqref{eq:one-sided} which is consistent with \cite[p. 504]{CT06}, in place of \cite[(3.7)]{CT06}, whose inconsistency with \cite[(3.10)]{CT06} is noted in \cite[Remark 3.9]{Hua20}.
\end{remark}

For $y\in(0,\delta)$, denote $Z_y:=\{F\leqslant y\}$. By \eqref{eq:one-sided} and \eqref{eq:band}, $Z_y$ is a compact region with smooth boundary, and since $\delta<\frac12$, it satisfies
\[
 \overline{B_{r+t_0}(p)}\subset Z_y\subset B_{r+t_0+a}(p)
 \qquad\text{and}\qquad
 \partial Z_y\subset U_{t_0}.
\]
Since $s/3<t_0$ and
$t_0+a\leqslant2s/3$, it satisfies $B_{r+s/3}(p)\subset Z_y\subset B_{r+2s/3}(p)$ and $\partial Z_y\subset A'$, and \cite[(3.6)]{CT06} and \eqref{eq:band} give $|\II|\leqslant|\operatorname{Hess}F|/|\nabla F|\leqslant C\ell^{-1}$ on $\partial Z_y$. We next choose the level $y_0$ by a second average. The coarea formula, \eqref{eq:band}, the inclusion $\{0<F<\delta\}\subset U_{t_0}$ and
\eqref{eq:choice-t} give
\begin{equation}\label{eq:cut-collar}
 \int_0^\delta\int_{\partial Z_y}w\,dA\,dy
 =\int_{\{0<F<\delta\}}w\,|\nabla F|\dvol
 \leqslant C\int_{U_{t_0}}\frac w\ell\dvol
 \leqslant\frac Cs\int_{A'}w\dvol .
\end{equation}
Therefore there exists $y_0\in(0,\delta)$ such that $Z:=Z_{y_0}$ satisfies \eqref{eq:weighted-cut}.
\end{proof}

\section{Estimating the boundary terms by the surrounding energy}\label{sec:scale}
Let $\varepsilon_0$ and $c_0$ be the constants of \cite[Theorem 0.8]{CT06}, and fix $L\geqslant4$ such that $c_0L^{-2}<\frac14$. Let $M$ be a complete noncompact Ricci-flat four-manifold, $p\in M$ and $0<s\leqslant r$. We take $a:=s/(512L)$ in \eqref{eq:radius}, put $A:=B_{r+s}(p)\setminus\overline{B_r(p)}$, and choose $Z$ with $\partial Z\subset A' \subset A$ by Lemma \ref{lem:weighted-cut} with weight $w=\ell/V_\ell$. 

\begin{proposition}\label{prop:cut}
The region $Z$ satisfies
\begin{equation}\label{eq:cut}
 b_0(\partial Z)+b_1(\partial Z)+J(\partial Z)
 \leqslant C\Bigl[\Bigl(\frac rs\Bigr)^3+\Bigl(\frac rs\Bigr)^{\frac34}\E(A)^{\frac34}\Bigr].
\end{equation}
\end{proposition}

\begin{proof}
Proposition \ref{prop:boundary-topology}, \eqref{eq:boundary} and \eqref{eq:weighted-cut} give the following control of the boundary terms
\begin{equation}\label{eq:cut-step}
 b_0(\partial Z)+b_1(\partial Z)+J(\partial Z)
 \leqslant C\int_{\partial Z}\frac{\ell}{V_\ell}\,dA
 \leqslant\frac Cs\int_{A'}\frac{\ell}{V_\ell}\dvol .
\end{equation}
For $x\in A'$, we have $B_{r+s}(p)\subset B_{4r}(x)$, so the Bishop-Gromov inequality at $p$ and then at $x$ gives
\begin{equation}\label{eq:volume-chain}
 \Vol A\leqslant\frac{Cs}r\Vol B_{4r}(x)\leqslant Csr^3\,\frac{V_\ell(x)}{\ell(x)^4},
 \qquad\text{so that}\qquad
 \frac{\ell(x)}{V_\ell(x)}\leqslant C\Bigl(\frac{sr^3}{\Vol A}\Bigr)^{\frac14}V_\ell(x)^{-\frac34}.
\end{equation}

$\bullet$ If $\ell(x)=a$, then \eqref{eq:volume-chain} and the Bishop-Gromov inequality give
\begin{equation}\label{eq:l=a}
    \frac{\ell(x)}{V_\ell(x)}\leqslant\frac{Csr^3}{a^3\Vol A}\leqslant\frac{Cr^3}{s^2\Vol A}.
\end{equation}
 
$\bullet$ If $\ell(x)<a$, then $\ell(x)=r_{|R|}(x)$ and $\int_{B_{L\ell(x)}(x)}|R|^2\dvol>\varepsilon_0$ by \cite[Theorem 0.8]{CT06}. By definition of the maximal function $M_f(x,t):=\sup_{0<\tau\leqslant t}\fint_{B_\tau(x)}|f|\dvol$, the Bishop-Gromov inequality, $\int_{B_{L\ell(x)}(x)}|R|^2\dvol>\varepsilon_0$ and $L\ell(x)<\frac s{512}$ give 
\begin{equation}\label{eq:V-m}
 \frac1{V_\ell(x)}\leqslant\frac{L^4}{\varepsilon_0}\fint_{B_{L\ell(x)}(x)}|R|^2\dvol
 \leqslant C\,M_{|R|^2}(x,\tfrac s{512}).
\end{equation}
Combining \eqref{eq:l=a} when $\ell =a$, with \eqref{eq:V-m} in \eqref{eq:volume-chain} when $\ell(x)<a$, gives in general:
\begin{equation}\label{eq:pointwise}
 \frac{\ell(x)}{V_\ell(x)}\leqslant C\Bigl[\frac{r^3}{s^2\Vol A}+\Bigl(\frac{sr^3}{\Vol A}\Bigr)^{\frac14}M_{|R|^2}(x,\tfrac s{512})^{\frac34}\Bigr].
\end{equation}
Finally, integrating \eqref{eq:pointwise} over $A'$ with $\Vol A'\leqslant\Vol A$, we obtain
\[
 \begin{aligned}
 \frac1s\int_{A'}\frac{\ell}{V_\ell}\dvol
 &\leqslant C\Bigl[\Bigl(\frac rs\Bigr)^3+\Bigl(\frac{r^3}{s^3\Vol A}\Bigr)^{\frac14}\int_{A'}M_{|R|^2}(x,\tfrac s{512})^{\frac34}\dvol(x)\Bigr]\\
 &\leqslant C\Bigl[\Bigl(\frac rs\Bigr)^3+\Bigl(\frac rs\Bigr)^{\frac34}\E(A)^{\frac34}\Bigr],
 \end{aligned}
\]
where we simplified the maximal function by applying \cite[(4.2)]{CT06} on $A$ with $W=A'$, $f=|R|^2$, radius $\frac s{512}$ and $\alpha=\frac34$. With \eqref{eq:cut-step}, this proves \eqref{eq:cut}.
\end{proof}

\begin{lemma}\label{prop:einstein}
Let $M$ be a complete noncompact Ricci-flat four-manifold, let $p\in M$ and $r>0$. If a neighborhood of $\overline{B_{2r}(p)}$ embeds smoothly into a four-manifold $N$ without boundary, then
\begin{equation}\label{eq:einstein-ball}
 \E(B_r(p))\leqslant b_2(N)+C.
\end{equation}
\end{lemma}

\begin{proof}
We may assume $b_2(N)<\E(B_r(p))$. Following the dyadic argument in \cite[proof of Theorem 1.26]{CT06}, we introduce for $i\geqslant0$,
\[
 r_i:=(2-2^{-i})r,\qquad
 A_i:=B_{r_{i+1}}(p)\setminus\overline{B_{r_i}(p)}
 \qquad\text{and}\qquad
 x_i:=2^{-3i}\bigl(\E(B_{r_i}(p))-b_2(N)\bigr).
\]
Since $0<x_i\leqslant2^{-3i}\E(B_{2r}(p))\to0$ as $i\to+\infty$, we can fix $i$ such that $x_i=\max_jx_j$, and then $x_{i+1}\leqslant x_i$ gives
\begin{equation}\label{eq:annulus-energy}
 2^{-3i}\E(A_i)=8x_{i+1}-x_i\leqslant7x_i .
\end{equation}
Choose the cut $Z$ of Proposition \ref{prop:cut} with $r_i$ and $2^{-i-1}r$ in place of $r$ and $s$, so that $A=A_i$. Since $B_{r_i}(p)\subset Z\subset B_{r_{i+1}}(p)$, the region $Z$ embeds in $N$. As $r_i/(r_{i+1}-r_i)\leqslant2^{i+2}$, the first inequality of \eqref{eq:reduction} with $\dim H^2_c(N)=b_2(N)$ by Poincaré duality \cite[Section 1.1.2]{Car07}, \eqref{eq:cut} and \eqref{eq:annulus-energy} give
\begin{equation}\label{eq:one-step}
 x_i\leqslant2^{-3i}\bigl(\E(Z)-b_2(N)\bigr)
 \leqslant C\bigl[1+\bigl(2^{-3i}\E(A_i)\bigr)^{\frac34}\bigr]
 \leqslant C\bigl(1+x_i^{\frac34}\bigr).
\end{equation}
If $x_i\geqslant1$, then \eqref{eq:one-step} gives $x_i\leqslant2Cx_i^{\frac34}$, that is, $x_i\leqslant(2C)^4$. Therefore $\E(B_r(p))-b_2(N)=x_0\leqslant x_i\leqslant\max\{1,(2C)^4\}$.
\end{proof}

\section{Proof of the main results}

\begin{proof}[Proof of Theorem \ref{thm:main}]\,

\emph{Energy bounds.} By \eqref{eq:einstein-ball} with $N=M$, we have
 $\E(M)=\sup_{r>0}\E(B_r(p))\leqslant b_2(M)+C.$
 
Assume now that $\E(M)<+\infty$ and fix compactly supported representatives of $m$ independent classes $([\omega_i])_{i\in\{1,\dots,m\}}$ in $H_c^2(M)$. Choose $r$ sufficiently large so that their supports lie in $B_r$ and $\E(M\setminus B_r)\leqslant1$. The region $Z_r$ of Proposition \ref{prop:cut} with $s=r$ then satisfies
\begin{equation}\label{eq:neighborhood-bound}
 B_r\subset\mathring{Z}_r\subset Z_r\subset B_{2r}
 \qquad\text{and}\qquad
 m\leqslant b_2(Z_r)\leqslant\E(M)+C.
\end{equation}
Indeed, extension by zero sends the restricted classes $[\omega_i|_{\mathring{Z}_r}]\in H^2_c(\mathring{Z}_r)$ to the original independent classes $[\omega_i]$ in $H_c^2(M)$, so the restricted classes are independent and $m\leqslant\dim H_c^2(\mathring{Z}_r)=b_2(Z_r)$ by Poincaré duality on $\mathring{Z}_r$ \cite[Section 1.1.2]{Car07}. 

(1) Since $\E(M)<+\infty$, applying \cite[Theorem 0.8]{CT06} to $B_{d(x,p)/2}(x)$ implies $|R(x)|\leqslant4c_0\,d(x,p)^{-2}$ for $d(x,p)$ large enough since  $\mathcal{E}(B_{d(x,p)/2}(x))$  tends to zero as $d(x,p)\to+\infty$. Together with collapsing, we are under the assumptions of \cite[Theorem 2.3]{CT06}, which gives the $N$-structure, and property iii) on \cite[p. 499]{CT06} gives orbits of extrinsic diameter $o(r)$. 

(2) For $\Sigma_r:=\partial Z_r$, $|R|\leqslant4c_0\,d(x,p)^{-2}$, \eqref{eq:weighted-cut} and $\Vol B_r(p)=o(r^4)$ give $J(\Sigma_r)\to0$. Since $\chi(Z_r)\in\mathbb Z$, \eqref{eq:gauss-bonnet} implies $\E(M)\in\mathbb Z$ and $\chi(Z_r)=\E(M)$ for every sufficiently large $r$.

(3) Assume $M$ is oriented. Then, as above, for $r$
large enough, a basis of $H_c^2(M)$ has representatives supported in $B_r(p)$, so extension by zero from $H_c^2(\mathring{Z}_r)$ is surjective and preserves the wedge pairing \cite[Section 1.1.2]{Car07}, hence $\tau(Z_r)=\tau(M)$.

Let $\eta$ denote the eta invariant of the odd signature operator on even forms. The boundary correction to \cite[Theorem 4.14]{APS75} is given by \cite[(2.6)]{DW07} and is
bounded by $CJ(\Sigma_r)\to0$. Since $|W^\pm|^2\leqslant |R|^2$ is integrable, the signature formula gives
\[
 \eta(\Sigma_r)\longrightarrow\eta_\infty
 =\frac1{12\pi^2}\int_M\bigl(|W^+|^2-|W^-|^2\bigr)\dvol-\tau(M).
\]
The last inequalities come from 
$|\eta_\infty+\tau(M)|\leqslant\frac{2}{3}\E(M)$, and
$|\tau(M)|\leqslant b_2(M)$.
\end{proof}

\begin{proof}[Proof of Corollary \ref{cor:exotic}]
If $X$ homeomorphic to $\R^4$ admits a complete Ricci-flat metric, then by \eqref{eq:einstein-ball} and \eqref{eq:neighborhood-bound}, every compact subset of $X$ is inside some $Z_r$ with $b_2(Z_r)\leqslant2C$. Hence $\gamma(X)\leqslant2C$ by \cite[Theorem 4.2]{Tay97}, which excludes the uncountably many smoothings of $\R^4$ with $\gamma=+\infty$ of \cite[Example 5.10]{Tay97}.
\end{proof}

\begin{proof}[Proof of Corollary \ref{cor:limits}]
For a limit $(M,g)$ and a fixed $r$, a neighborhood of $\overline{B_{2r}(p)}$ embeds in $N_j$ for all large $j$. Applying \eqref{eq:einstein-ball}, taking the lower limit in $j$, and then letting $r\to+\infty$ shows the bound for smooth metrics. The orbifold version is deduced from the orbifold Chern-Gauss-Bonnet formula (see e.g.
\cite{Nak90} for the adaptation to ALE ends): each singular point $q$
contributes $1-|\Gamma_q|^{-1}$ to the corrected energy.
\end{proof}

\end{document}